\documentclass{daj}
\usepackage{amsmath}
\usepackage{amssymb}
\usepackage{bbm}

\numberwithin{equation}{section}

\newtheorem{thm}{Theorem}[section]
\newtheorem{cor}[thm]{Corollary}
\newtheorem{lem}[thm]{Lemma}
\newtheorem{prop}[thm]{Proposition}

\newcommand{\ii}[2][\delta']{\left\langle #2 \right
\rangle^{#1}}
\def\conv{\mathrm{conv}}

\def\d{\mathtt{d}}

\dajAUTHORdetails{%
  title = {Floating and Illumination Bodies for Polytopes: Duality Results}, 
  author = {Olaf Mordhorst and Elisabeth M. Werner},
  plaintextauthor = {Olaf Mordhorst, Elisabeth M. Werner},
  runningtitle = {Floating and Illumination Bodies for Polytopes: Duality Results}, 
  runningauthor = {Olaf Mordhorst and Elisabeth M. Werner},
  copyrightauthor = {O. Mordhorst, E. M. Werner},
  keywords = {floating body, illumination body},
}   

\dajEDITORdetails{%
   year={2019},
   number={11},
   received={11 September 2017},   
   revised={7 December 2018},    
   published={18 June 2019},  
   doi={10.19086/da.8973},       
}   

\begin{document}

\begin{frontmatter}[classification=text]

\title{Floating and Illumination Bodies for Polytopes: Duality Results} 

\author[OlMo]{Olaf Mordhorst\thanks{Partially supported by the German Academic Exchange Service and DFG project BE 2484/5-2}}
\author[ElWe]{Elisabeth M. Werner\thanks{Partially supported by  NSF grants DMS-1504701, 811146}}

\begin{abstract}
We consider the question how well a floating body can be approximated by the polar of the illumination body of the polar. We establish precise convergence results in the case of centrally symmetric polytopes. 
This leads to a new affine invariant which is related to the cone measure of the  polytope.
\end{abstract}
\end{frontmatter}

\section{Introduction}

Floating bodies are a fundamental notion in convex geometry. Early notions  of floating bodies are motivated by the physical description of floating objects. The first systematic study of floating bodies appeared 1822 in a work by C. Dupin \cite{Dupin:1822} on naval engineering. Blaschke, in dimensions $2$ and $3$,  and later Leichtweiss in higher dimensions, used the floating body in the study
of affine differential geometry,  in particular affine surface area (see \cite{Blaschke}, \cite{Leichtweiss:1986}).
Affine surface area is among the most powerful tools in convexity.  It  is widely used, for instance  in  approximation of convex bodies by polytopes, e.g.,  \cite{Boeroetzky2000/2, Reitzner:2010, Schuett:2003}
and affine differential geometry, e.g.,
 \cite{Andrews:1999, Ivaki:2013a, Stancu:2003, TW2}.  
\par
Dupin's floating body may not exist and related to that,  originally affine surface area was only defined for sufficiently smooth bodies. 
Sch\"utt and Werner \cite{SchuettWerner1990} and independently B{\'a}r{\'a}ny and Larman \cite{Barany:1988}  introduced the convex floating body which, 
in contrast to Dupin's  original definition,  always exists  and coincides with Dupin's floating
body if  the latter exists.
This,  in turn,  allowed to define affine surface area for all convex bodies as carried out  in \cite{SchuettWerner1990}.
\par
The illumination body was introduced much later  in \cite{Werner94} as a further tool to study affine properties of convex bodies.
It was pointed out in \cite{Werner06} that the convex floating body and the illumination body are dual notions in the sense that the polar of the floating body and the illumination body of the polar should be close both,  in a conceptual and geometric sense. In \cite{MordhorstWerner1} the duality relation is studied in detail for \(C_+^2\)-bodies, i.e. bodies with twice differentiable boundary and everywhere positive curvature,  and \(\ell_p^n\)-balls. The paper provides asymptotic sharp estimates how well the polar of the floating body
can be approximated by an illumination body of the polar.  A limiting procedure leads to  a new affine invariant that is  different from the affine surface area. It is related  to the cone measure of the convex body.  These measures play a central role in many aspects of convex geometry,  e.g.,  \cite{BLYZ, CLM, Naor, PaourisWerner2012}.
\vskip 3mm
The purpose of this paper is to make the duality relation between floating body and illumination body precise  when the convex body is a polytope $P$. 
It was shown by Sch\"utt  \cite{Schuett91} that the limit of the (appropriately normalized) volume difference of a polytope $P$ and its floating body leads to a quantity related to the combinatorial structure of the
polytope, namely the flags  of $P$ (see section \ref{Combinatorics}). Likewise,  as shown in \cite{Werner96}, the limit of the (appropriately normalized) volume difference of a simplex  and its illumination  body is related to the combinatorics of the boundary. 
Now, as in the smooth case \cite{MordhorstWerner1}, a limit procedure leads to a new  affine invariant that is not related to the  combinatorial structure of the boundary of the polytope, but, as in the smooth case,  to cone measures.
\vskip 2mm
The techniques in \cite{MordhorstWerner1} rely on comparing ``extremal" directions, i.e., 
 directions where the boundary of the convex body and its floating body, and  the illumination body of its polar, differ the most and the least. 
The techniques used  in \cite{MordhorstWerner1} employ tools from differential geometry which is possible due to the  \(C_+^2\) smoothness assumptions.
Such tools are no longer available in our present setting  of polytopes  and we 
have to use a completely different approach.
\vskip 2mm
It would be interesting to have results for more general classes of convex bodies other than  polytopes and the ones investigated in \cite{MordhorstWerner1}. A major obstruction is that in general it is hard to compute the polar body and even harder to  
compute  the floating body if we do not have smoothness assumptions or are in the case of polytopes.
The smooth case \cite{MordhorstWerner1} and the polytope case seem to be  the  extremal cases. 
Indeed, ellipsoids are  \(C_+^2\)-bodies where equality of the polar of the floating body and the illumination body of the polar can be achieved. 
And we show here that polytopes display the  worst behavior  
one can expect. Furthermore, the limit is not continuous with respect to the polytope involved since it depends only on the local structure of the boundary (see section \ref{Approximation}). 
\vskip 2mm
The  paper is organized as follows. In the next section we present the main theorem and some  consequences. In Section \ref{Tools+Lemmas} we give 
the necessary background material and several lemmas needed for the proof of the main theorem.
In section \ref{Combinatorics} we discuss properties of the new affine invariant. 
We show, with an example, that it  is not continuous with respect to  the Hausdorff distance. We also show that 
for this invariant the combinatorial  structure of the polytope is less relevant. The relation to the cone measures is the dominant feature.
In the final section we address questions of  approximation of the floating body by the polar of the illumination of the polar. We show that our convergence results are of pointwise nature and we derive a uniform upper bound for general convex bodies.

\section{Main theorem and consequences}

We work in a similar framework as in \cite{MordhorstWerner1}. We first recall the  notions and definitions that we will need. 
Let \(K\) be a \(n\)-dimensional convex body and \(\delta \geq 0\). The convex floating body $K_\delta$ of $K$ was introduced by Sch\"utt and Werner \cite{SchuettWerner1990} and independently by Barany and Larman  \cite{Barany:1988} as
\begin{equation}\label{floatingbody}
K_{\delta}= \bigcap\limits_{|K \cap H^{-}|_n\leq \delta |K|_n} H^+ , 
\end{equation}
where \(H\) is a hyperplane and \(H^+, H^-\) are the corresponding closed half-spaces and $|K|_n$ is the volume of $K$.
\par
\noindent
We denote by \(\conv[A,B]\) the convex hulls of two sets \(A,B\subseteq \mathbb{R}^n\). If \(B=\{x\}\) we simply write \(\conv[A,x]\)
 for the  convex hull of $A$ and the vector $x$. 
\par
\noindent 
The  illumination body $K^\delta$ of  $K$ was defined by  E. Werner  \cite{Werner94} as  
\begin{equation}\label{illuminationbody}
K^{\delta}=\{x\in\mathbb{R}^n: |\mathrm{conv}[K,x]|_n\leq (1+\delta)|K|_n\}.
\end{equation}
\par
\noindent
If \(0\) is in the interior of \(K\),  the polar \(K^{\circ}\) of \(K\) is  given by
\begin{equation} \label{polar}
K^\circ=\{y \in \mathbb{R}^n: \langle x, y\rangle \leq 1,    \text{ for all } x \in K\}.
\end{equation}
It is a simple fact of convex geometry that for a hyperplane \(H\) with corresponding halfspaces \(H^+, H^{-}\) there is a corresponding point \(x_H\in\mathbb{R}^n\) such that
\[
\left(K\cap H^{+}\right)^{\circ}=\conv[K^{\circ},x_H], 
\]
whenever \(0\) is in the interior of \(K\cap H^+\).  As noted in \cite{Werner06},  this polarity relation gives rise to the idea that cutting off with hyperplanes sets of a certain  volume of a convex body and including points such that their convex hull with the convex body has a certain volume, should be dual operations, 
\begin{align}
\left(K_{\delta}\right)^{\circ}\approx \left(K^{\circ}\right)^{\delta'}\label{gleich}.
\end{align}
In the same paper it is pointed out that equality cannot be achieved in general since the floating body is always strictly convex and the illumination body of a polytope is always a polytope.  
\par
\noindent
As in \cite{MordhorstWerner1} we like to measure how ``close" these two bodies are in the polytope case.  A further outcome  of such a  study shows how well the floating body of a polytope can be approximated by a polytope, namely the polar of an illumination of the polar.
For \(x\in\mathbb{R}^n\backslash\{0\}\) and \(K\) a convex body with \(0\in K\) we denote by \(r_K(x)=\sup\{\lambda\geq 0: \lambda x\in K\}\) the radial function of \(K\). To measure how close two centrally symmetric convex bodies  \(S_1,S_2\) are, we use the distance 
\begin{equation}\label{distance}
\d(S_1,S_2)=\sup_{u\in\mathbb{S}^{n-1}}\max\left[\frac{r_{S_1}(u)}{r_{S_2}(u)},\frac{r_{S_2}(u)}{r_{S_1}(u)}\right]=\inf \left\{a\geq 1: \frac{1}{a}S_1\subseteq S_2 \subseteq a S_1\right\}.
\end{equation}
It is worthwhile to mention that \(\log \d(\cdot,\cdot)\) is a metric  on the space of convex bodies in $\mathbb{R}^n$ which induces the same topology as the Hausdorff distance. 
\par
\noindent
For   a centrally symmetric convex body  \(S\) and \(0<\delta<\frac{1}{2}\), we put \(\ii[\delta]{S}= \left(\left(S^{\circ}\right)^{\delta}\right)^{\circ}\). We then define
\begin{equation} \label{FloatIlluMeas}
\d_S(\delta)=\inf_{\delta'>0}\d\left(S_{\delta},\ii{S}\right) . 
\end{equation}
Please note   that \(\d_{L(S)}(\delta)=\d_S(\delta)\) for every linear invertible map \(L\).
Observe also that $\d_{B_2^n}(\delta) = 1$ . 
\vskip 3mm
\noindent
One of the main theorems in \cite{MordhorstWerner1} states that for origin symmetric convex bodies $C$  in $\mathbb{R}^n$   that are 
$C^2_+$, i.e. the Gauss curvature $\kappa (x)$ exists for every $x \in \partial C$ and is strictly positive, the relation (\ref{gleich}) can be made precise  in terms of the cone measures
of $C$ and $C^\circ$. 
\par
For a Borel set $A \subset \partial C$,  the cone measure $M_C$ of $A$  is  defined as $M_C(A) = | \conv[0, A]|_n$. The density function of  $M_C$ is $m_{C} (x) = \frac{1}{n}  \langle x, N(x)\rangle$ (see \cite{PaourisWerner2012}),    and we write  $n_{C} (x) = \frac{1}{n |C|_n}  \langle x, N(x)\rangle $ for  the density of the normalized cone measure $\mathbb{P}_C$ of $C$ (again  see e.g.,  \cite{PaourisWerner2012}). 
This means  that, e.g., \cite{PaourisWerner2012}, 
$$
d M_C (x) = m_{C} (x) d   \mu_C(x) \  \   \text {and } \  \    d \mathbb{P}_C (x) = n_{C} (x) d   \mu_C(x).
$$
\par
\noindent
Denote by $N_C : \partial C \rightarrow S^{n-1}$, $x \rightarrow N(x)$ the Gauss map of $C$, see e.g.,  \cite{SchneiderBook}. Then, 
similarly, $m_{C^\circ} (x) = \frac{1}{n}  \frac{\kappa_C(x) } {\langle x, N(x)\rangle ^n} $ is the density function of the ``cone measure"  $M_{C^\circ}$ of $C^\circ$. 
For a Borel set $A \subset \partial C$,   $M_{C^\circ} (A) = | \conv[0, N_{C^\circ}^{-1} (N_C(A))]|_n$ and $n_{C^\circ} (x) = \frac{1}{n |C^\circ|_n}  \frac{\kappa_C(x) } {\langle x, N(x)\rangle^n} $ is the density of the normalized cone measure $\mathbb{P}_{C^\circ}$ of $C^\circ$, see, e.g., \cite{PaourisWerner2012}. When $C$ is $C^2_+$,  this formula holds for 
all $x \in \partial C$.
Thus
$$
d   M_{C^\circ} (x) = m_{C^\circ} (x) d \mu_C(x) \  \   \text {and } \  \    d \mathbb{P}_{C^\circ} (x) = n_{C^\circ} (x) d \mu_C(x).
$$
As observed in \cite{MordhorstWerner1},  we then have for a centrally symmetric $C^2_+$ convex body $C$ that
\begin{eqnarray*} \label{Theorem:smooth}
 \lim _{\delta \rightarrow \infty} \frac{\d_C(\delta) -1}{\delta^{\frac{2}{n+1}} }
 = c_n \  \left(|C|_n|C^{\circ}|_n\right)^{\frac{1}{n+1}}  \  \left[ \max_{x\in\partial C} \left(\frac{n_{C^\circ} (x)}{n_{C} (x) } \right)^{\frac{1}{n+1}}  - \min_{x\in\partial C} \left(\frac{n_{C^\circ} (x)}{n_{C} (x) } \right)^{\frac{1}{n+1}}  \right],
 \end{eqnarray*}
where $c_n = \frac{1}{2} \ \left(\frac{n+1}{|B_2^{n-1}|_{n-1}}\right)^{\frac{2}{n+1}}$.
\vskip 3mm
\noindent
In the case of a polytope the  (discrete) densities $n_{P}$ and  $n_{P^\circ}$ of the  normalized cone measures    
can be expressed as follows. 
Let $\xi$ be an extreme point of $P$. Let $F_\xi$ be the facet of $P^\circ$ that has $\xi$ as an outer normal.  The (discrete) density   of the  normalized cone measure of $P^\circ$ is
\begin{equation} \label{cone-poly-polar}
n_{P^\circ} (\xi)  = \frac{1}{n \  |P^\circ|_n}  \  \frac{1}{\|\xi\|} \  |F_\xi|_{n-1},
\end{equation}
where $\| \cdot\|$ denotes the standard Euclidean norm on $\mathbb{R}^n$.
 Let \(s(F_{\xi})\) be the \((n-1)\)-dimensional Santal\'{o} point (see, e.g., \cite{Gardner1995, SchneiderBook}) of \(F_{\xi}\) and \((F_{\xi}-s(F_{\xi}))^{\circ}\) be the polar of \((F_{\xi}-s(F_{\xi}))\) with respect to the \((n-1)\)-dimensional subspace in which \(F_{\xi}-s(F_{\xi})\) lies. We put
\begin{equation} \label{cone-poly}
n_{P}(\xi)  = \frac{1}{n \  |P|_n}  \  \|\xi\| \  |(F_{\xi}-s(F_{\xi}))^{\circ}|_{n-1} .
\end{equation} 
Let $C_{\xi}$ be the cone with base $F_{\xi}$ and apex at the origin and let 
$$
C_{\xi}^* = \{ y \in \mathbb{R}^n: \langle x, y \rangle \geq 0,   \text{ for all } x \in C_{\xi}\}
$$ 
be the cone dual to $C_{\xi}$. Then $|(F_{\xi}-s(F_{\xi}))^{\circ}|_{n-1}$ is the $(n-1)$-dimensional volume of the base of
\[
Z=C_{\xi}^*\cap\left\{x\in\mathbb{R}^n:\left\langle x,\frac{\xi}{\|\xi\|}\right\rangle\leq \xi\right\}
\] 
and thus \(\frac{1}{n}\|\xi\| |(F_{\xi}-s(F_{\xi}))^{\circ}|_{n-1}\) is the \(n\)-dimensional volume of the finite cone \(Z\). The expression $n_{P}(\xi)$ is the ratio of this volume and the volume of $P$. We see $n_{P}(\xi)$ as a cone measure  associated to \(\xi\) since this volume is the cone measure of the set of all points of \(Z\) with outer normal \(\frac{\xi}{\|\xi\|}\).
\vskip 3mm
\noindent
Our main theorem can be  expressed in terms of  $n_{P}(\xi)$ and $n_{P^\circ}(\xi)$ and reads as follows.
\begin{thm}  \label{MT}

Let \(P\subseteq \mathbb{R}^n\) be a centrally symmetric polytope. Then
\[
\lim_{\delta\rightarrow 0}\frac{\d_P(\delta)-1}{\delta^{1/n}}= \min_{c \geq 0}  \left[  \max_{\xi\in\mathrm{ext}(P)} \left(\frac{n_{P^\circ} (\xi) - c\  n_{P} (\xi)^\frac{1}{n}}{n_{P} (\xi)^\frac{1}{n} \  n_{P^\circ} (\xi)} , \ \frac{c}{\min_{\eta\in\mathrm{ext}(P)} n_{P^\circ} (\eta)} \right)  \right].
\] 
\end{thm}
\vskip 2mm
Recall that for $1 \leq p < \infty$,  the \(\ell_p^n\)- unit balls are defined as $B^n_p=\{x \in \mathbb{R}^n: \left(\sum_{i=1}^ n |x_i|^p\right)^\frac{1}{p}  \leq 1\} $.
The subsequent corollary about  the cube $B^n_\infty = \{x \in \mathbb{R}^n: \max_{1 \leq i \leq n} |x_i| \leq 1\}$ and the crosspolytope $B^n_1=\{x \in \mathbb{R}^n: \sum_{i=1}^ n |x_i| \leq 1\} $,  is an immediate consequence of Theorem \ref{MT}.
\vskip 2mm
\begin{cor} \label{MC}
\[
\lim_{\delta\rightarrow 0}\frac{\d_{B_{\infty}^n}(\delta)-1}{\delta^{1/n}}=\frac{\sqrt[n]{n!}}{n} \   \   \text{and} \  \  \lim_{\delta\rightarrow 0}\frac{\d_{B_1^n}(\delta)-1}{\delta^{1/n}}=\frac{2^{1/n}}{2} . 
\]
\end{cor}

\section{Tools and Lemmas }\label{Tools+Lemmas}
Let \(P\subseteq \mathbb{R}^n\) be a centrally symmetric polytope. In \cite{MeyerReisner} it was shown that for centrally symmetric convex bodies 
Dupin's floating body exists and coincides with the convex floating body. This means that every support hyperplane of \(P_{\delta}\) cuts off the volume \(\delta|P|_n\) from \(P\). We use this fact throughout the remainder of the paper.  We denote by $\mathrm{ext}(P)$ the set of extreme points of $P$. Note that this set coincides with the set of vertices of $P$.  For \(\xi\in \mathrm{ext}(P) \), let  \(F_1,\dots,F_k\) be the \((n-1)\)-dimensional facets of \(P\) such that \(\xi\in F_i\). Then there are \(y_1,\dots,y_k\in\mathbb{R}^n\) such that 
for \(1\leq i\leq k\),
\[
F_i\subseteq H_i:=\{x\in\mathbb{R}^n: \langle x, y_i\rangle =1\}.
\]
Observe that \(y_1,\dots,y_k\) are vertices of \(P^{\circ}\) and that \(F_{\xi}:=\conv[y_1,\dots,y_k]\) is a facet of \(P^{\circ}\). 
Let \(s(F_{\xi})\) be the \((n-1)\)-dimensional Santal\'{o} point of \(F_{\xi}\) and \((F_{\xi}-s(F_{\xi}))^{\circ}\) be the polar of \((F_{\xi}-s(F_{\xi}))\) with respect to the \((n-1)\)-dimensional subspace in which \(F_{\xi}-s(F_{\xi})\) lies (see (\ref{cone-poly-polar}) and (\ref{cone-poly})). 
\vskip 2mm
For $\delta >0$, let $P_\delta$ be the floating body of $P$. Let \(\xi\in \mathrm{ext}(P) \). 
 We denote by \(\xi_{\delta}\) the unique point in the intersection of \(\partial P_{\delta}\) with the line segment \([0,\xi]\) and by \(\ii[\delta]{x}\) the unique point in the intersection of \(\partial \ii[\delta]{P}\) with \([0,\xi]\). We denote by \(\xi^{\delta}\) the unique point such that \(\xi\) is the unique point in the intersection of \(\partial P^{\delta}\) with \([0,\xi^{\delta}]\).\\

The next lemma  provides a formula for 
\(\frac{\|\xi_{\delta}\|}{\|\xi\|}\)  if \(\delta>0\) is sufficiently small.

\begin{lem}\label{FloatingPointPolytope} Let \(P\subseteq \mathbb{R}^n\) be a centrally symmetric polytope. Then there is  \(\delta_{0}>0\) such that for every \(0\leq \delta\leq \delta_0\) and every vertex \(\xi\in\partial P\) we have
\[
\frac{\|\xi_{\delta}\|}{\|\xi\|}=1-\left(\frac{n|P|_{n}}{|(F_{\xi}-s(F_{\xi}))^{\circ}|_{n-1}\|\xi\|}\right)^{1/n}\delta^{1/n}.
\]
\end{lem}
\textit{Proof.} Let \(e_i\in\mathbb{R}^n\) be the vector with \(i\)-th entry \(1\) and the other entries are \(0\). We first consider the case that \(\xi=e_n\) and \(s(F_{\xi})=e_n\). For \(v\in\mathbb{R}^n\backslash \{0\}\) we denote by \(v^{\perp}=\{w\in\mathbb{R}^n:\langle v,w\rangle =0\}\) the orthogonal complement of \(v\). We show that 
\[
e_n^{\perp}\cap\bigcap_{i=1}^k \{x\in\mathbb{R}^n:\langle x, y_i\rangle\leq 1\} 
\]
is an \((n-1)\)-dimensional convex body with centroid in the origin. For self similarity reasons the \((n-1)\)-dimensional centroid of 
\[
(\alpha e_n+ e_n^{\perp})\cap\bigcap_{i=1}^k \{x\in\mathbb{R}^n:\langle x, y_i\rangle\leq 1\} 
\]
is \(\alpha e_n\) for every \(\alpha<1\). Let \(\bar{y}_i\in\mathbb{R}^{n-1}\) be such that \((\bar{y}_i,1)=y_i\), \(1\leq i\leq k\). Put  
\[
F=\{\bar{y}\in\mathbb{R}^{n-1}: (\bar{y},1)\in F_{e_n}\}=\conv[\bar{y}_1,\dots,\bar{y}_k]\subseteq \mathbb{R}^{n-1}
\]
and 
\[
B=F^{\circ}=\bigcap_{i=1}^k \{\bar{x}\in\mathbb{R}^{n-1}: \langle\bar{x},\bar{y}_i\rangle\leq 1\}
\]
where the polar is taken in  \(\mathbb{R}^{n-1}\). Then \(s(F)=0\) and
\[
\bigcap_{i=1}^k \{x\in\mathbb{R}^n:\langle x, y_i\rangle\leq 1\} 
=\{(\lambda\bar{x},1-\lambda): \lambda \geq 0,  \  \bar{x}\in B\}.
\]
It is a well-known fact (see \cite{Santalo}) that for a convex body \(C\) we have the identities 
\[
g\left((C-s(C))^{\circ}\right)=0=s\left((C-g(C))^{\circ}\right).
\]
It follows immediately that the centroid of 
\[
e_n^{\perp}\cap \bigcap_{i=1}^k \{x\in\mathbb{R}^n:\langle x, y_i\rangle\leq 1\}
\]
lies in the origin. The volume of the cone with base \(B\times\{0\}\) and apex \(\xi=e_n\) is given by \(|B|_{n-1}/n\). Let \(0\leq \Delta\leq 1\). Then the volume of the cone with base
\[
((1-\Delta)e_n+e_n^{\perp})\cap\bigcap_{i=1}^k\{x\in\mathbb{R}^n:\langle x,y_i\rangle\leq 1\}
\]
and apex \(e_n\) is given by \(\frac{\Delta^n}{n}|B|_{n-1}\). 
There is  \(\Delta_0>0\) such that for every \(0\leq \Delta\leq \Delta_0\) the point \(e_n\) is the only vertex of \(P\) contained in the half-space
\[
\{x\in\mathbb{R}^n:x_n\geq 1-\Delta\}.
\]
Hence, the above described cone is given by
\[
\{x\in\mathbb{R}^n:x_n\geq 1-\Delta\}\cap P.
\]
Let \(\delta>0\) and choose \(\Delta\) such that
\[
\frac{1}{n}|B|_{n-1}\Delta^{n}=\delta|P|_n,
\]
or, equivalently, 
\[
\Delta=\left(\frac{n|P|_n}{|B|_{n-1}}\right)^{\frac{1}{n}}\delta^{1/n}.
\]
Choose \(\delta_0>0\) sufficiently small such that  for every \(0\leq \delta\leq \delta_0\) the value of \(\Delta\) is smaller than or equal to \(\Delta_0\). 
It was shown in \cite{MeyerReisner} that for centrally symmetric convex bodies, the floating body coincides with the convex floating body.
Thus, since \(P\) is centrally symmetric, 
the floating body of $P$ coincides with the convex floating body, and therefore  the hyperplane 
\(
\{x\in\mathbb{R}^n:x_n\geq 1-\Delta\}
\)
touches \(P_{\delta}\) at the centroid of
\[
\{x\in\mathbb{R}^n:x_n\geq 1-\Delta\}\cap P.
\]
This centroid is given by
\[
(1-\Delta)e_n=\left(1-\left(\frac{n|P|_n}{|(B^{\circ}-s(B^{\circ}))^{\circ}|_{n-1}}\right)^{\frac{1}{n}}\delta^{1/n}\right)e_n=\left(1-\left(\frac{n|P|_n}{|(F-s(F))^{\circ}|_{n-1}\|e_n\|}\right)^{\frac{1}{n}}\delta^{1/n}\right)e_n.
\]
\vskip 2mm
\noindent
For a general vertex \(\xi\) and general \(s(F_{\xi})\) note first that \(\langle \xi,s(F_{\xi})\rangle=1\) and thus, \(\xi\not\in s(F_{\xi})^{\perp}\). Let \(L\in\mathbb{R}^{n\times n}\) be a matrix with last row \(s(F_{x_i})\) and the other rows are a basis of \(\xi^{\perp}\). Let \(L^{-t}\) the inverse of the transpose. Since \(\langle\xi,s(F_{\xi})\rangle=1\) it follows that \(L\) is a full rank matrix with \(L(\xi)=e_n\) and \(L^{-t}(s(F_{\xi}))=e_n\). Then, \(LP\) is a centrally symmetric polytope with vertex \(L(\xi)=e_n\) and \(s(F_{L(\xi)})=s(F_{e_n})=e_n\). Note that \(\frac{\|\xi_{\delta}\|}{\|\xi\|}=\frac{\|(L(\xi))_{\delta}\|}{\|L(\xi)\|}\). The lemma follows from

\begin{eqnarray*}
\frac{n|LP|_n}{|(F_{L\xi}-s(F_{L\xi}))^{\circ}|_{n-1}\|L\xi\|}
&=& \  \frac{n|\det(L)|\cdot|P|_n}{|(L^{-t}(F_{\xi}-s(F_{\xi})))^{\circ}|_{n-1}\|L\xi\|}\\
&= &\frac{n|\det(L)|\cdot|P|_n}{\left(|\det(L^{-t})|\cdot\|L\frac{\xi}{\|\xi\|}\|\right)^{-1}   |(F_{\xi}-s(F_{\xi})))^{\circ}|_{n-1}\|L\xi\|}\\
&= &\frac{n|P|_n}{|\left(F_{\xi}-s(F_{\xi})\right)^{\circ}|_{n-1}\|\xi\|}\quad.
\end{eqnarray*}
The second equality follows from the fact that for every \((n-1)\)-dimensional vector space \(V\) with normal \(u\), every linear invertible map \(S:\mathbb{R}^n\rightarrow\mathbb{R}^n\) and every Borel set \(A\subseteq V\),  we have
\(|S(A)|_{n-1}=|\det(S)|\cdot \|S^{-t}(u)\|\cdot|A|_{n-1}\).
\hfill \(\Box\)
\vskip 3mm
\noindent
For a vertex \(\xi\in P\),  \(\ii[\delta]{\xi}\) is the unique point in the intersection of  \( \partial \left(\ii[\delta]{P}\right) \)  and the line segment \([0,\xi]\). 

\vskip 2mm
\begin{lem}\label{IlluminationPointPolytope} Let \(P\) be a centrally symmetric polytope. Then there is a \(\delta_0>0\) such that for every \(0\leq\delta\leq \delta_0\) and every extreme point \(\xi\in \mathrm{ext}(P)\) we have
\[
\frac{\|\ii[\delta]{\xi}\|}{\|\xi\|}=\left(1+\frac{n|P^{\circ}|_n\|\xi\|}{|F_{\xi}|_{n-1}}\delta\right)^{-1}.
\]
\end{lem}
\textit{Proof.} We show that 
\[
\left\{y\in\mathbb{R}^{n}:\langle y, \xi\rangle=1+\frac{n|P^{\circ}|_n\|\xi\|}{|F_{\xi}|_{n-1}}\delta\right\}
\]
is a support hyperplane of \((P^{\circ})^{\delta}\). The lemma then follows immediately. \\
Let \(z\in F_{\xi}\) and \(\Delta\geq 0\). The volume of the cone with base \(F_{\xi}\) and apex \(z+\Delta\frac{\xi}{\|\xi\|}\) is \(\frac{1}{n}|F_{\xi}|_{n-1}\Delta\). There is a \(\Delta_0>0\) and an \(\eta>0\) such that
\[
F_{\xi}^{\eta}=\{z\in F_{\xi}: \mathrm{dist}(z,\partial F_{\xi})\geq \eta\}
\]
has non-empty relative interior and such that for every \(z\in F_{\xi}^{\eta}\) and every \(0\leq \Delta\leq \Delta_0\) we have
\[
\conv\left[P^{\circ}, z+\Delta\frac{\xi}{\|\xi\|}\right]=P^{\circ}\cup\conv\left[F_{\xi},z+\Delta\frac{\xi}{\|\xi\|}\right].
\]
Let \(\delta_0>0\) be such that \(\frac{n|P^{\circ}|_n}{|F_{\xi}|_{n-1}}\delta\leq \Delta_0\) for every \(0\leq \delta\leq \delta_0\). 
It is obvious  that for every \(z\in F_{\xi}^{\eta}\) the vector 
\[
z+\frac{n|P^{\circ}|_n}{|F_{\xi}|_{n-1}}\delta\frac{\xi}{\|\xi\|}
\]
lies on the boundary of \((P^{\circ})^{\delta}\). Since \(F_{\xi}\) is contained in the hyperplane \(\{y\in\mathbb{R}^n: \langle y, \xi\rangle = 1\}\),  it follows that 
\[
\left\{y\in\mathbb{R}^{n}:\langle y, \xi\rangle=1+\frac{n|P^{\circ}|_n\|\xi\|}{|F_{\xi}|_{n-1}}\delta\right\}
\]
is a support hyperplane of \((P^{\circ})^{\delta}\).

\hfill \(\Box\)

\begin{lem} Let \(P\subseteq \mathbb{R}^n\) be a centrally symmetric polytope. Then there is a \(\delta_0>0\) such that for every \(0\leq \delta\leq \delta_0\)
\begin{align*}
\conv[\{\ii[\delta]{\xi}:\xi\in \mathrm{ext}(P)\}]\subseteq \ii[\delta]{P}
\subseteq \conv\left[\left\{\ii[\delta]{\xi}:\xi\in \mathrm{ext}(P)\right\}\cup\left\{\frac{1}{2}(\xi+\xi'):\xi,\xi'\in\mathrm{ext}(P):\xi\neq \xi'\right\}\right].
\end{align*}
\end{lem}
\textit{Proof.} The first inclusion is obvious. Choose \(\delta_0>0\) as in Lemma \ref{IlluminationPointPolytope}. \\
The second inclusion will follow from the fact that  \((P^{\circ})^{\delta}\supseteq C(\delta)\), where
\[
C(\delta)= \bigcap_{\xi\in\mathrm{ext}(P)}\left\{y\in\mathbb{R}^{n}:\langle y, \xi\rangle=1+\frac{n|P^{\circ}|_n\|\xi\|}{|F_{\xi}|_{n-1}}\delta\right\}\cap
\bigcap_{\xi,\xi'\in\mathrm{ext}(P),\xi\neq \xi'}\{y\in\mathbb{R}^{n}:\frac{1}{2}\langle\xi+\xi',y\rangle\leq 1\}.
\]
Let \(y_0\in C(\delta)\backslash P^{\circ}\). It follows that there is a \(\xi\in \mathrm{ext}(P)\) with \(\langle y_0, \xi\rangle\geq 1\). Let \(\xi'\in\mathrm{ext}(P)\backslash\{\xi\}\). Since \(\langle y_0, \frac{1}{2}\left(\xi+\xi'\right)\rangle\leq 1\) it follows that
\[
\langle y_0, \xi'\rangle = \langle y_0, \xi+\xi'\rangle -\langle y_0, \xi\rangle
\leq 2-1=1. 
\]
We obtain
\[
\conv[P^{\circ}, y_0]=P^{\circ}\cup \conv[F_{\xi}, y_0] .
\]
Since \(y_0 \in \left\{y\in\mathbb{R}^{n}:\langle y, \xi\rangle=1+\frac{n|P^{\circ}|_n\|\xi\|}{|F_{\xi}|_{n-1}}\delta\right\}\),  we deduce
\[
\left|\conv[F_{\xi}, y_0]\right|_{n}=\frac{1}{n}\left(\left\langle y_0,\frac{\xi}{\|\xi\|}\right\rangle-\frac{1}{\|\xi\|}\right)|F_{\xi}|_{n-1}\leq \frac{1}{n}\cdot\frac{n|P^{\circ}|_n}{|F_{\xi}|_{n-1}}\delta\cdot|F_{\xi}|_{n-1}=\delta|P^{\circ}|,
\]
which means that  \(y_0\in (P^{\circ})^{\delta}\).
\hfill \(\Box\)
\vskip 2mm
\noindent
\begin{cor}\label{BoundsPolarIlluminationBody} Let \(P\subseteq \mathbb{R}^n\) be a centrally symmetric polytope. Then there is a \(\delta_0>0\) such that for every \(0\leq \delta\leq \delta_0\)
\begin{align*}
&\conv[\{\ii[\delta]{\xi}:\xi\in \mathrm{ext}(P)\}]\subseteq \ii[\delta]{P}\\
\subseteq & 
\conv\left[\left\{\ii[\delta]{\xi}:\xi\in \mathrm{ext}(P)\right\}\cup\left\{\frac{1}{2}(\xi+\xi'):\xi,\xi'\in\mathrm{ext}(P):\xi\neq \xi'\text{ and }\frac{1}{2}(\xi+\xi')\in\partial P\right\}\right].
\end{align*}

\end{cor}
\textit{Proof.} We only need to prove 
\[
\ii[\delta]{P}
\subseteq 
\conv\left[\left\{\ii[\delta]{\xi}:\xi\in \mathrm{ext}(P)\right\}\cup\left\{\frac{1}{2}(\xi+\xi'):\xi,\xi'\in\mathrm{ext}(P):\xi\neq \xi'\text{ and }\frac{1}{2}(\xi+\xi')\in\partial P\right\}\right].
\]
Consider the set 
\[
\left\{\frac{1}{2}(\xi+\xi'):\xi,\xi'\in\mathrm{ext}(P):\xi\neq \xi'\text{ and }\frac{1}{2}(\xi+\xi')\in\mathrm{int} (P)\right\}.
\]
This set is finite and \(\lim\limits_{\delta\rightarrow 0}\ii[\delta]{\xi}=\xi\) for every \(\xi\in\mathrm{ext}(P)\). It follows that there is a \(\delta_0>0\) such that for every \(0\leq\delta\leq \delta_0\) the following holds
\[
\left\{\frac{1}{2}(\xi+\xi'):\xi,\xi'\in\mathrm{ext}(P):\xi\neq \xi'\text{ and }\frac{1}{2}(\xi+\xi')\in\mathrm{int} (P)\right\}\subseteq \conv\left[\left\{\ii[\delta]{\xi}:\xi\in\mathrm{ext}(P)\right\}\right]
\]
which yields the claim of the corollary.
\hfill \(\Box\)
\vskip 2mm
\noindent
\begin{lem}\label{UpperBoundFloatingBody} Let \(P\subseteq \mathbb{R}^n\) be a centrally symmetric polytope. Then there is a function \(t:[0,\frac{1}{2}]\rightarrow \mathbb{R}\) with \(\lim_{\delta\rightarrow 0}t(\delta)=0\) such that 
\[
P_{\delta}\subseteq(1-\Lambda\delta(1-t(\delta)))P, 
\]
where \(\Lambda=\min_{\zeta\in\mathrm{ext}(P^{\circ})} \frac{|P|_n\|\zeta\|} {|F_{\zeta}|_{n-1}}\).
\end{lem}
\textit{Proof.} Let $\delta >0$ be given. Let $\zeta\in\mathrm{ext}(P^{\circ})$. We choose  \(\Delta=\Delta(\zeta,\delta) \) such that 
$$\left|P\cap \left\{x\in\mathbb{R}^n: \left\langle x,\frac{\zeta}{\|\zeta\|}\right\rangle\geq \frac{1}{\|\zeta\|}-\Delta\right\}\right|_n = \delta |P|_n .
$$
For $\delta>0$ and hence
 \(\Delta=\Delta(\zeta,\delta) \geq 0\)  sufficiently small,  the volume of 
\(P\cap\{x\in\mathbb{R}^n: \langle x,\frac{\zeta}{\|\zeta\|}\rangle\geq \frac{1}{\|\zeta\|}-\Delta\}\) is up to an error given by \(\Delta|F_{\zeta}|_{n-1}\), i.e. there is a function \(T_{\zeta}\) with \(\lim_{\Delta\rightarrow 0}T_{\zeta}(\Delta)=0\) such that
\[
\left|\left\{x\in\mathbb{R}^n: \left\langle x,\frac{\zeta}{\|\zeta\|}\right\rangle\geq \frac{1}{\|\zeta\|}-\Delta  \right\} \cap P\right|_n=\Delta|F_{\zeta}|_{n-1}(1+T_{\zeta}(\Delta)).
\]
Hence, for every \(\zeta\in\mathrm{ext}(P^{\circ})\),  there is a  function \(t_{\zeta}\) with \(\lim_{\delta\rightarrow 0}t_{\zeta}(\delta)=0\) such that  
\[
P_{\delta}\subseteq \left\{x\in\mathbb{R}^n:\langle\zeta, x\rangle\leq 1-\frac{|P|_n \|\zeta\|}{|F_{\zeta}|_{n-1}}\delta(1- t_{\zeta}(\delta))\right\}.
\]
Let $t(\delta) = \max_{\zeta\in\mathrm{ext}(P^{\circ})}t_\zeta (\delta)$ and $\Lambda=\min_{\zeta\in\mathrm{ext}(P^{\circ})} \frac{|P|_n \|\zeta\|} {|F_{\zeta}|_{n-1}}$. Then 
\[
P_{\delta}\subseteq \bigcap_{\zeta\in\mathrm{ext}(P^{\circ})} \left\{x\in\mathbb{R}^n:\langle\zeta, x\rangle\leq 1-\Lambda\delta(1-t(\delta))\right\}=(1-\Lambda\delta(1-t(\delta)))P.
\]
\hfill \(\Box\)

\begin{lem}\label{UnextremalPointsFlBody} Let \(P\subseteq \mathbb{R}^n\) be a centrally symmetric polytope and \(x\in\partial P\backslash \mathrm{ext}(P)\). Then there exists  \(\delta_0>0\) and  \(k>0\) such that for every \(0\leq \delta\leq \delta_0\) we have
\[
\frac{\|x_{\delta}\|}{\|x\|}\geq  1-k\  \delta^{\frac{1}{n-1}}.
\]
\end{lem}
\textit{Proof.} Since \(x\) is not an extreme point of $P$, there are points \(x_1,x_2\in\partial P\) with \(x_1\neq x\neq x_2\) and such that \(x=\frac{1}{2}(x_1+x_2)\). By a linear transformation of \(P\) we may assume without loss of generality that \(x=e_2\),  \(x_1=e_2-e_1\) and  \(x_2=e_2+e_1\). There is an \(0<\varepsilon<1\)
such that \([-\varepsilon,\varepsilon]\times \{0\}\times [-\varepsilon, \varepsilon]^{n-2}\subseteq P\). It follows that the centrally symmetric convex body
\[
S=\conv\left[e_2\pm\varepsilon e_1, -e_2\pm\varepsilon e_1, [-\varepsilon,\varepsilon]\times \{0\}\times [-\varepsilon, \varepsilon]^{n-2}\right]=[-\varepsilon,\varepsilon]\times \conv\left[\pm e_2,\{0\}\times\{0\}\times [-\varepsilon,\varepsilon]^{n-2}\right]
\]
is contained in \(P\). Put \(\tilde{\delta}=\delta\frac{|P|_n}{|S|_n}\). We compute \((e_2)_{\tilde{\delta}}\) with respect to \(S_{\tilde{\delta}}\). Let $0 \leq \Delta <1$. A simple computation shows that 
\[
|S\cap \{x\in\mathbb{R}^n: x_2\geq 1-\Delta \}|_n=\frac{1}{n}(2\varepsilon\Delta)^{n-1}.
\]
The \((n-1)\)-dimensional centroid of the $(n-1)$-dimensional set 
\(
S\cap \{x\in\mathbb{R}^n: x_2=1-\Delta \}
\)
lies on the line \(\mathbb{R}e_2\). Since $S$ is symmetric,  the convex floating and the floating body of Dupin coincide \cite{MeyerReisner} and   it follows that for \(\delta<\frac{1}{2}\), 
\[
\left(1-\frac{(n|S|_n)^{\frac{1}{n-1}}}{2\varepsilon}\tilde{\delta}^{\frac{1}{n-1}}\right)e_2\in \partial\left (S_{\tilde{\delta}}\right).
\]
 Since \(S\subseteq P\),  there exists  \(\delta_0>0\) and \( k>0\) such that for every \(0\leq\delta\leq\delta_0\), 
\[
\frac{\|x_{\delta}\|}{\|x\|}\geq 1-k \  \delta^{\frac{1}{n-1}}
\]
where \(x_{\delta}\) is taken with respect to \(P_{\delta}\).
\hfill \(\Box\)

\section{Proof of Theorem \ref{MT}  and Corollary \ref{MC} } \label{Beweise}

We recall the quantities that are relevant for our main theorem. For $\xi \in \mathrm{ext}(P)$, we put 
\begin{equation} \label{alpha}
\alpha_{\xi}=\left(\frac{n|P|_n}{|(F_{\xi}-s(F_{\xi}))^{\circ}|_{n-1}\|\xi\|}\right)^{1/n}, 
\end{equation}
\begin{equation} \label{beta}
\beta_{\xi}=\frac{n|P^{\circ}|_n\|\xi\|}{|F_{\xi}|_{n-1}}\  \   \text{and}  \  \  \beta=\max_{\xi\in\mathrm{ext}(P)}\beta_{\xi} .
\end{equation}
For $ c\geq 0$, we set 
\begin{equation} \label{gamma}
G_c(P)=\max_{\xi\in\mathrm{ext}(P)}[a_{\xi}-c\beta_{\xi},c\beta] \  \   \text{and}  \  \ G(P)=\min_{c\geq 0}G_c(P).
\end{equation}
\vskip 3mm
\noindent
Then  Theorem \ref{MT} reads.
\vskip 2mm
\noindent
{\bf Theorem \ref{MT}} 
{\em  Let \(P\subseteq \mathbb{R}^n\) be a centrally symmetric polytope. Then
\[
\lim_{\delta\rightarrow 0}\frac{\d_P(\delta)-1}{\delta^{1/n}}= G(P).
\] 
}
\vskip 3mm
\noindent
We split  the proof of the theorem  and  show separately   the upper and lower bound.

\subsection{Upper bound}
We prove the following proposition.
\begin{prop} Let \(P\subseteq \mathbb{R}^n\) be a centrally symmetric polytope. Then 
\[
\limsup_{\delta\rightarrow 0}\frac{\d_P(\delta)-1}{\delta^{1/n}}\leq G(P).
\]
\end{prop}
\textit{Proof.} Let \(c_0\geq 0\) be such that \(G(P)=G_{c_0}(P)\) and put \(\delta'=c_0 \ \delta^{1/n}\). 
By Lemma \ref{IlluminationPointPolytope}, Lemma \ref{BoundsPolarIlluminationBody} and Lemma \ref{UpperBoundFloatingBody}, a sufficient condition for \(P_{\delta}\subseteq a\ii{P}\) is that
\[
1-\Lambda\delta(1-t(\delta))\leq a(1+\beta_{\xi}\delta')^{-1}, 
\]
for every \(\xi\in\mathrm{ext}(P)\). Hence, 
\[
a\geq (1-\Lambda\delta(1-t(\delta)))\max_{\xi\in\mathrm{ext}(P)}(1+\beta_{\xi}\delta')=(1-\Lambda\delta(1-t(\delta)))(1+\beta c_o\delta^{1/n}).
\]
By  Lemma \ref{FloatingPointPolytope}, Lemma \ref{IlluminationPointPolytope} and  Corollary \ref{BoundsPolarIlluminationBody}, a sufficient condition for \(\ii{P}\subseteq a P_{\delta}\) is that
\[
(1+\beta_{\xi}\delta')^{-1}\leq a(1-\alpha_{\xi}\delta^{1/n})
\]
for every \(\xi\in\mathrm{ext}(P)\) and that 
\[
\left\|\frac{1}{2}(\xi+\xi')\right\|\leq a\left\|\left(\frac{1}{2}(\xi+\xi')\right)_{\delta}\right\|
\]
for \(\xi,\xi'\in\mathrm{ext}(P)\), \(\xi\neq\xi'\) such that \(\frac{1}{2}(\xi+\xi')\in\partial P\). From the first condition we derive that
\[
a\geq \frac{1}{(1-\alpha_{\xi}\delta^{1/n})(1+\beta_{\xi}\delta')}
\]
for every \(\xi\in\mathrm{ext}(P)\). 
By Lemma \ref{UnextremalPointsFlBody} there is a constant \(k>0\) and  \(\delta_0>0\) such that for every \(0\leq\delta\leq \delta_0\) we have
\[
\left\|\left(\frac{1}{2}(\xi+\xi')\right)_{\delta}\right\|\geq \left(1-k \ \delta^{\frac{1}{n-1}}\right) \  \left\|\frac{1}{2}(\xi+\xi')\right\|
\] 
and we may assume that \(k\) and \(\delta_0\) are taken uniformly with respect to all pairs \((\xi,\xi')\). Hence, for \(\delta\leq\delta_0\) we have the condition that 
\[
a\geq \frac{1}{1-k \  \delta^{\frac{1}{n-1}}}.
\]
We check that all three conditions are met if one takes \(a=1+G(P) \ \delta^{\frac{1}{n}}\ (1+o(1))\).
The condition
\[
a\geq \frac{1}{1-k\delta^{\frac{1}{n-1}}}
\]
is obvious since \(1+G(P)\delta^{1/n}\geq (1-k\delta^{1/{n-1}})^{-1}\) for sufficiently small \(\delta>0\). The condition
\[
a\geq\frac{1}{(1-\alpha_{\xi}\delta)(1+\beta_{\xi}\delta')}
\]
is true since
\begin{align*}
&\frac{1}{(1-\alpha_{\xi}\delta^{1/n})(1+\beta_{\xi}\delta')}
=\frac{1}{(1-\alpha_{\xi}\delta^{1/n})(1+\beta_{\xi}c_0\delta^{1/n})}
= 1+(\alpha_{\xi}-c_0\beta)\delta^{1/n}+o(\delta^{1/n})\\
\leq & 1+G_{c_0}(P)\delta^{1/n}+o(\delta^{1/n})
\leq 1+G(P)\delta^{1/n}+o(\delta^{1/n}).
\end{align*}
Finally, the condition
\[
a\geq (1-\Lambda\delta(1-t(\delta)))(1+\beta c_o\delta^{1/n})
\]
is true since
\[
(1-\Lambda\delta(1-t(\delta)))(1+\beta c_o\delta^{1/n})\leq 1+c_0 \beta\delta^{1/n}\leq 1+G_{c_0}(P)\delta^{1/n}\leq 1+G(P)\delta^{1/n}+o(\delta^{1/n}).
\]

\hfill \(\Box\)

\subsection{Lower Bound}
We prove the following proposition.
\begin{prop} Let \(P\subseteq \mathbb{R}^n\) be a centrally symmetric polytope. Then 
\[
\liminf_{\delta\rightarrow 0}\frac{\d_P(\delta)-1}{\delta^{1/n}}\geq G(P).
\]
\end{prop}
\textit{Proof.} Let \(c_0 \geq  0\) such that \(G(P)=G_{c_0}(P)\) and let \(\xi_1,\xi_2\in\mathrm{ext}(P)\) be such that
\[
\beta_{\xi_1}=\max_{\zeta\in\mathrm{ext}(P)}\beta_{\zeta}\quad \text{ and }\quad
\alpha_{\xi_2}-c_0\beta_{\xi_2}=\max_{\zeta\in\mathrm{ext}(P)}[\alpha_{\zeta}-c_0\beta_{\zeta}].
\]
We obtain that \(c_0\beta_{\xi_1}=\alpha_{\xi_2}-c_o\beta_{\xi_2}\) and  therefore that  \(c_0=\frac{\alpha_{\xi_2}}{\beta_{\xi_1}+\beta_{\xi_2}}\) and \(G(P)=\frac{\alpha_{\xi_2}\beta_{\xi_1}}{\beta_{\xi_1}+\beta_{\xi_2}}\).
A necessary condition for \(\ii{P}\subseteq aP_{\delta}\) is that \(\|\ii{\xi_2}\|\leq a\|(\xi_2)_{\delta}\|\), or, equivalently, also using Lemmas \ref{FloatingPointPolytope} and \ref{IlluminationPointPolytope}, 
\[
a\geq \frac{\|\ii{\xi_2}\|}{\|(\xi_2)_{\delta}\|} = (1-\alpha_{\xi_2}\delta^{1/n})^{-1}(1+\beta_{\xi_2}\delta')^{-1}.
\]
By Lemma \ref{BoundsPolarIlluminationBody},  there is  \(\delta_0\) such that for every \(0\leq\delta\leq\delta_0\) we have
\[
\ii{P}
\subseteq \conv\left[\left\{\ii{\xi}:\xi\in \mathrm{ext}(P)\right\}\cup\left\{\frac{1}{2}(\xi+\xi'):\xi,\xi'\in\mathrm{ext}(P):\xi\neq \xi'\right\}\right]=:P(\delta').
\]
If \(\delta'_0>0\) is chosen sufficiently small,  \(\ii{\xi_1}\) is an extreme point of \(P(\delta')\).  Then there exists  \(\varepsilon>0\) and a hyperplane \(H^y=\{x\in\mathbb{R}^n:\langle x, y\rangle= 1\}\)
such that \(\left\langle\ii{\xi_1},y\right\rangle>1+\varepsilon\) and such that all  other extreme points of \(P(\delta')\) lie in \(\{x\in\mathbb{R}^n:\langle x,y \rangle\leq 1\}\) for every \(0\leq\delta'\leq\delta'_0\). Hence, 
\[
\ii{P}\cap\{x\in\mathbb{R}^n:\langle x,y\rangle\geq 1\}
\subseteq \conv\left[P\cap H^y, \ii{\xi_1}\right].
\] 
Let \(z\in(\partial P)\cap H^y\).  Then \(\lambda \ii{\xi_1}+(1-\lambda)z\not\in\mathrm{int}\left[\ii{P}\right]\),  for every \(\lambda\in[0,1]\). Fix \(\lambda\in [0,1]\) and 
put \(v=\lambda\xi_1+(1-\lambda)z\in\partial P\). Let \(t,\mu\in (0,1)\) be such that \(tv=\mu\ii{\xi_1}+(1-\mu)z\).  Then
\begin{equation}\label{gleichung:v}
t\|v\|\geq \|\ii{v}\|.
\end{equation}
We determine \(t\). By Lemma \ref{IlluminationPointPolytope} and as $\xi$ and $\ii{\xi_1}$, we know that 
\[
\ii{\xi_1}=(1+\beta_{\xi_1}\delta')^{-1}\xi_1
\]
if \(\delta'_0>0\) is chosen sufficiently small. This means that  \(t\) and \(\mu\) satisfy the equation
\[
t(\lambda\xi_1+(1-\lambda)z)=\mu(1+\beta_{\xi_1}\delta')^{-1}\xi_1+(1-\mu)z.
\]
Since \(\xi\) and \(z\) are linearly independent,  \(t\) and \(\mu\) satisfy the system of linear equations
\begin{align*}
\mathrm{I.}\quad t\lambda-\mu(1+\beta_{\xi_1}\delta')^{-1}=0\qquad\quad\mathrm{II.}\quad t(1-\lambda)+\mu=1.
\end{align*}
It follows that \(t=(1+\lambda\beta_{\xi_1}\delta')^{-1}\). Since \(v\) is not an extreme point  of $P$, it follows from Lemma \ref{UnextremalPointsFlBody} that there is a \(k_v\geq 0\) such that
\[
\frac{\|v_\delta\|}{\|v\|}\geq 1-k_v \ \delta^{\frac{1}{n-1}}.
\]
By this and (\ref{gleichung:v}), a necessary condition for \(a\ii{P}\supseteq P_{\delta}\) is that 
\(
a(1+\lambda\beta_{\xi_1}\delta')^{-1}\geq 1-k_v \ \delta^{\frac{1}{n-1}}
\), i.e.,  \(a\geq (1+\lambda\beta_{\xi_1}\delta')( 1-k_v\ \delta^{\frac{1}{n-1}})\). Assume that \(\delta'\geq \frac{\alpha_{\xi_2}}{\lambda \beta_{\xi_1}+\beta_{\xi_2}}\delta^{1/n}\)       
then we get
\[
a\geq (1+\lambda\beta_{\xi_1}\delta')( 1-k_v\ \delta^{\frac{1}{n-1}})\geq
1+\frac{\alpha_{\xi_2}\beta_{\xi_1}\lambda}{\lambda\beta_{\xi_1}+\beta_{\xi_2}}\delta^{1/n}+o(\delta^{1/n}).
\] 
The assumption \(\delta'\leq \frac{\alpha_{\xi_2}}{\lambda \beta_{\xi_1}+\beta_{\xi_2}}\delta^{1/n}\) together with the necessary condition 
\(a\geq (1-\alpha_{\xi_2}\delta^{1/n})^{-1}(1+\beta_{\xi_2}\delta')^{-1}\) also yields
\[
a\geq 1+\frac{\alpha_{\xi_2}\beta_{\xi_1}\lambda}{\lambda\beta_{\xi_1}+\beta_{\xi_2}}\delta^{1/n}+o(\delta^{1/n}).
\]
Thus,
\[
\liminf_{\delta\rightarrow 0}\frac{\d_P(\delta)-1}{\delta^{1/n}}\geq \frac{\alpha_{\xi_2}\beta_{\xi_1}\lambda}{\lambda\beta_{\xi_1}+\beta_{\xi_2}}.
\]

Letting \(\lambda\rightarrow 1\),  we get the desired result.

\hfill \(\Box\)

\vskip 3mm
\subsection{ Proof of Corollary \ref{MC}}

We first treat the case of the cube.

\begin{cor}
\[
\lim_{\delta\rightarrow 0}\frac{\d_{B_{\infty}^n}(\delta)-1}{\delta^{1/n}}=\frac{\sqrt[n]{n!}}{n}
\]
\end{cor}
\textit{Proof.} By symmetry,  \(\alpha_{\xi}\) and \(\beta_{\xi}\) have the same value for  all the extreme points of \(B_{\infty}^n\). Take \(\xi=(1,\dots,1)\). Then
\(\|\xi\|=\sqrt{n}\), \(|B_{\infty}^n|_n=2^n\), \(|B_1^n|_n=\frac{2^n}{n!}\) and 
\[
|F_{\xi}|_{n-1}=|\conv[e_1,\dots,e_n]|_{n-1}=\frac{\sqrt{n}}{(n-1)!}.
\]
It is well known that the volume product  $\left|S_{n-1}\right|_{n-1} \left|  (S_{n-1})^\circ\right|_{n-1} $ of the \((n-1)\)-dimensional simplex is \(\frac{n^n}{((n-1)!)^2}\). Hence,
as $F_{\xi}$ is an $(n-1)$-dimensional regular simplex, 
\[
|(F_{\xi}-s(F_{\xi}))^{\circ}|_{n-1}=\frac{1}{|F_\xi|_{n-1}}\cdot \frac{n^n}{((n-1)!)^2}=\frac{n^n}{\sqrt{n}(n-1)!}.
\]
Therefore,
\[
\alpha_{\xi}=\left(\frac{n2^n}{\frac{n^n}{\sqrt{n}(n-1)!}\sqrt{n}}\right)^{1/n}=2\frac{\sqrt[n]{n!}}{n} \  \  \text{and}  \  \   \beta_\xi=2^n .
\]
The minimum over all $c\geq 0$ of \(\max[\alpha_{\xi}-c\beta_{\xi},c\beta_{\xi}]\) is attained for \(c=\frac{\alpha_{\xi}}{2\beta_{\xi}}\). Thus 
\[
G(B_{\infty}^n)=\frac{\alpha_{\xi}}{2}= \frac{\sqrt[n]{n!}}{n} , 
\]
which completes the proof.

\hfill \(\Box\)
\vskip 3mm
Now we show the statement of Corollary \ref{MC} in the case of the crosspolytope.

\begin{cor}
\[
\lim_{\delta\rightarrow 0}\frac{\d_{B_1^n}(\delta)-1}{\delta^{1/n}}=\frac{2^{1/n}}{2}.
\]
\end{cor}
\textit{Proof.} As in the previous example,  all \(\alpha_{\xi}\) and all \(\beta_{\xi}\)
are equal and \(G(B_1^n)=\frac{\alpha_{\xi}}{2}\). Take \(\xi=e_n\). Then \(|B_1^n|_n=\frac{2^n}{n!}\), \(\|\xi\|=1\) and  \(F_{\xi}=\conv[e_n+\sum_{i=1}^{n-1}\theta_i e_i:\theta\in\{-1,1\}^{n-1}]=e_n+B_{\infty}^{n-1}\). It follows that
\[
|(F_{\xi}-s(F_{\xi}))^{\circ}|_{n-1}=|B_1^{n-1}|_{n-1}=\frac{2^{n-1}}{(n-1)!}.
\]
We obtain 
\[
\alpha_{\xi}=\left(\frac{n\frac{2^n}{n!}}{\frac{2^{n-1}}{(n-1)!}\cdot 1}\right)^{1/n}
=2^{1/n}.
\]
\hfill \(\Box\)

\section{The combinatorial structure of $\d_P$}\label{Combinatorics}

In \cite{Schuett91},  it  was proved  that the following relation holds for  all polytopes \(P\subseteq \mathbb{R}^n\), 
\[
\lim_{\delta\rightarrow 0}\frac{|P|_n-|P_{\delta}|_n}{\delta\mathrm{ln}(\delta)^{n-1}}= \frac{\operatorname{fl}_{n}(P)}{n!n^{n-1}},
\]
where $\operatorname{fl}_{n}(P)$ denotes the number of flags of $P$. A flag of $P$ is an 
$(n+1)$-tuple $(F_{0},\dots,F_{n})$ such that $F_{i}$ is an $i$-dimensional face of
$P$ and $F_{0}\subset F_{1}\subset\cdots\subset F_{n}$.
\par
This theorem suggests that also  \(\d_P\),  and hence $G(P)$,  might only depend on the combinatorial structure of \(P\). The fact that \(\d_P\) is invariant under affine transformations of \(P\) supports this conjecture. However, this is not the case,   as is illustrated by the following \(2\)-dimensional example. 
\vskip 2mm
\noindent
For \(\varepsilon\in (0,1)\), we consider the hexagon
\[
P(\varepsilon)=\conv\left[\pm e_2, \pm \sqrt{1-\varepsilon^2}e_1\pm \varepsilon e_2\right].
\]
We show that \(\d_{P(\varepsilon)}\) changes for different values of \(\varepsilon\). We compute the \(2\)-dimensional volume of \(P(\varepsilon)\). The hexagon is,  up to a nullset,  the disjoint union of the two congruent trapezoids
\begin{align*}
T_1&=\conv[-e_2, \sqrt{1-\varepsilon^2}e_1-\varepsilon e_2,\sqrt{1-\varepsilon^2}e_1+\varepsilon e_2,e_2] 
\end{align*}
and
\begin{align*}
T_2&=\conv[e_2, -\sqrt{1-\varepsilon^2}e_1+\varepsilon e_2,-\sqrt{1-\varepsilon^2}e_1-\varepsilon e_2,-e_2].
\end{align*}
The trapezoid \(T_1\) has the two parallel sides \(S_1=\conv[-e_2,e_2]\) and \(S_2=\conv[ \sqrt{1-\varepsilon^2}e_1-\varepsilon e_2,\sqrt{1-\varepsilon^2}e_1+\varepsilon e_2]\) and the height of \(T_1\) with respect to \(S_1, S_2\) is given by \(\sqrt{1-\varepsilon^2}\). Hence,
\[
|T_1|_2=\frac{|S_1|_1+|S_2|_1}{2}\cdot \sqrt{1-\varepsilon^2}=\frac{2+2\varepsilon}{2}\cdot\sqrt{1-\varepsilon^2}=(1+\varepsilon)\cdot\sqrt{1-\varepsilon^2},
\]
and we conclude that \(|P(\varepsilon)|_2=2\cdot|T_1|_2=2(1+\varepsilon)\cdot\sqrt{1-\varepsilon^2}\). \\
\par
\noindent
We compute the vertices of the polar of \(P(\varepsilon)\). One vertex is given as the solution of the equations
\[
y_2=1  \  \  \text{and} \  \   \sqrt{1-\varepsilon^2}y_1+\varepsilon y_2=1, 
\]
which yields \((y_1,y_2)=\left(\frac{1-\varepsilon}{\sqrt{1-\varepsilon^2}},1\right)\). Another vertex is given as the solution of the equations
\[
\sqrt{1-\varepsilon^2}y_1+\varepsilon y_2=1   \  \  \text{and} \  \    \sqrt{1-\varepsilon^2}y_1+\varepsilon y_2=1, 
\]
which yields \((y_1,y_2)=\left(\frac{1}{\sqrt{1-\varepsilon^2}}, 0\right)\). By symmetry,  the six vertices of \(P^{\circ}\) are given by 
\[
\left\{\pm \frac{1}{\sqrt{1-\varepsilon^2}}e_1, \pm\frac{1-\varepsilon}{\sqrt{1-\varepsilon^2}}e_1\pm e_2\right\}.
\]
Since \(P(\varepsilon)^{\circ}\) is the union of two trapezoids,  computations similar to the case of \(P(\varepsilon)\) yield that the \(2\)-dimensional volume of \(P(\varepsilon)^{\circ}\) is given by
\[
|P(\varepsilon)^{\circ}|_2=\frac{4-2\varepsilon}{\sqrt{1-\varepsilon^2}}.
\]
If \(\xi=\pm e_2\),  we get that  \(|F_{\xi}|_1=2\cdot\frac{1-\varepsilon}{\sqrt{1-\varepsilon^2}}\) and 
\[
|(F_{\xi}-s(F_{\xi}))^{\circ}|_1=2\cdot \left(\frac{|F_{\xi}|_1}{2}\right)^{-1}=2\frac{\sqrt{1-\varepsilon^2}}{1-\varepsilon}.
\]
Hence,
\[
\alpha_1:=\alpha_{\xi}=\left(\frac{2\cdot 2(1+\varepsilon)\cdot \sqrt{1-\varepsilon^2}}{2\cdot\frac{\sqrt{1-\varepsilon^2}}{1-\varepsilon}\cdot 1}\right)^{1/2}=\sqrt{2}\cdot\sqrt{1-\varepsilon^2}
\]
and
\[
\beta_1:=\beta_{\xi}=\frac{2\cdot \frac{4-2\varepsilon}{\sqrt{1-\varepsilon^2}}\cdot 1}{2\cdot\frac{1-\varepsilon}{\sqrt{1-\varepsilon^2}}}=\frac{4-2\varepsilon}{1-\varepsilon}.
\]
If \(\xi=\pm\sqrt{1-\varepsilon^2}e_1\pm\varepsilon e_2\) then 
\[
|F_{\xi}|_1=\left\|\frac{1-\varepsilon}{\sqrt{1-\varepsilon^2}}e_1+e_2-\frac{1}{\sqrt{1-\varepsilon^2}}e_1\right\|=\frac{1}{\sqrt{1-\varepsilon^2}}
\]
and
\[
|(F_{\xi}-s(F_{\xi}))^{\circ}|_1=2\cdot\left(\frac{1}{2\cdot\sqrt{1-\varepsilon^2}}\right)^{-1}=4\cdot\sqrt{1-\varepsilon^2}.
\]
Hence,
\[
\alpha_2:=\alpha_{\xi}=\left(\frac{2\cdot 2(1+\varepsilon)\cdot\sqrt{1-\varepsilon^2}}{4\cdot \sqrt{1-\varepsilon^2}}\right)^{1/2}=(1+\varepsilon)^{1/2}
\]
and
\[
\beta_2:=\beta_{\xi}=\frac{2\cdot \frac{4-2\varepsilon}{\sqrt{1-\varepsilon^2}}}{\frac{1}{\sqrt{1-\varepsilon^2}}}=8-4\varepsilon.
\]
\par
\noindent
We compute \(G(P(\varepsilon))\). If \(0 < \varepsilon <\frac{1}{2}\),  then \(8-4\varepsilon>\frac{4-2\varepsilon}{1-\varepsilon}\) and therefore, \(\beta=\max_{\xi\in\mathrm{ext}(P_{\varepsilon})}\beta_{\xi}=8-4\varepsilon\). Moreover,  for \(0 < \varepsilon <\frac{1}{2}\),
\(\alpha_1>\alpha_2\) and thus \(\alpha_1-c\cdot\beta_1\geq \alpha_2-c\cdot\beta_2\),  for every \(c\geq 0\). This yields
\[
G_c(P(\varepsilon))=\max\left[c(8-4\varepsilon), \sqrt{2}\cdot \sqrt{1-\varepsilon^2}-c\cdot\frac{4-2\varepsilon}{1-\varepsilon}\right]
\]
and \(G_c(P(\varepsilon))\) is minimized by
\[
c_0=\frac{(1+\varepsilon)^{1/2} \ (1-\varepsilon)^{3/2}}{\sqrt{2}\  (2-\varepsilon)\ (3-2\varepsilon)}.
\]
It follows that
\[
G(P(\varepsilon))=G_{c_0}(P(\varepsilon))= 2\ \sqrt{2}\cdot\frac{(1+\varepsilon)^{1/2}(1-\varepsilon)^{3/2}}{3 - 2  \varepsilon}.
\]
This means that,   if \(\varepsilon>0\) is sufficiently small,  \(G(P(\varepsilon))\)  and hence \(\d_{P(\varepsilon)}\), changes for different values of \(\varepsilon\).
\par
\noindent
Moreover, this example shows that the affine invariant $G(P(\varepsilon))$ is not continuous with respect to the Hausdorff distance,  since \(P(\varepsilon)\) converges to \(B_1^2\) as \(\varepsilon\) goes to \(0\) but 
\[
\lim_{\varepsilon\rightarrow 0}2\ \sqrt{2}\cdot\frac{(1+\varepsilon)^{1/2}(1-\varepsilon)^{3/2}}{3 - 2  \varepsilon}
=\frac{2\ \sqrt{2}}{3}\neq \frac{\sqrt{2}}{2}=G(B_1^2).
\]

\section{Approximation results for the floating body and open questions}\label{Approximation}

The parameter \(\d_S\) measures the best approximation of the floating body by the polar of an illumination body of the polar. We establish a uniform bound for this quantity, independent of the convex body.
\begin{prop}\label{ApproximationConvexBodyFloatingBody} Let \(S\subseteq \mathbb{R}^n\) be a centrally symmetric convex body. Then there exist constants \(G_n, \delta_n\) only depending on the dimension such that 
\[
\d(S_{\delta}, S)\leq 1 + G_n \delta^{1/n}
\]
for \(\delta\in [0,\delta_n]\). 
\end{prop}
In particular, the proposition yields
\[
\d_S(\delta)=\inf_{\delta'\geq 0} \d\left(S_\delta,\ii{S}\right)
\leq \d\left(S_{\delta},\ii[0]{S}\right)=\d(S_{\delta}, S)\leq 1+G_n\delta^{1/n}.
\]
Thus,  \(\delta^{1/n}\) is already the worst order of convergence we can get in general.  The polytopal case shows that we cannot hope for any better uniform rate of convergence. Proposition \ref{ApproximationConvexBodyFloatingBody} does not involve the floating body any more. We address the question if we get a better uniform bound if we involve the illumination body.  At best, how does the optimal \(G_n\) look like such that \(\d_S(\delta)\leq 1+ G_n\delta^{1/n}+o(\delta^{1/n})\) where the error term \(o(\delta^{1/n})\) does only depend on the dimension?
It would be interesting to know about the maximizers if they exist. Furthermore, it would be interesting to know something about the best uniform bound on subclasses like polytopes, \(C^1\)-bodies or \(C^2\)-bodies.
\vskip 2mm
\noindent
\textit{Proof of Proposition \ref{ApproximationConvexBodyFloatingBody}.} The quantity \(\d(S_{\delta},S)\) is invariant with respect to linear transformations and we may therefore assume that \(S\) is in John position, i.e.,  in particular,
\[
B_2^n\subseteq S\subseteq \sqrt{n}B_2^n.
\]
Let \(\xi\in\partial S\).  Then \(1\leq \|\xi\|\leq \sqrt{n}\). Put \(B_2^{\xi^{\perp}}=B_2^n\cap \xi^{\perp}\). By central symmetry of  \(S\), the double cone  $C_{\xi}=\conv[\xi,B_2^{\xi^{\perp}},-\xi]$ is contained in $S$, 
\[
C_{\xi}=\conv[\xi,B_2^{\xi^{\perp}},-\xi]\subseteq S.
\]
 Put 
\[
\Delta=\left(\frac{n|S|_n}{\|\xi\||B_2^{n-1}|_{n-1}}\right)^{1/n}\delta^{1/n}.
\]
The halfspace
\[
\{x\in\mathbb{R}^n: \langle x, \xi\rangle\geq 1-\Delta\}
\]
cuts off exactly volume \(\delta|S|_n\) from \(C_{\xi}\). Similar to the proof of Lemma \ref{UnextremalPointsFlBody} we get that
\[
\frac{\|\xi_{\delta}\|}{\|\xi\|}\geq 1-\left(\frac{n|S|_n}{\|\xi\||B_2^{n-1}|_{n-1}}\right)^{1/n}\delta^{1/n}.
\]
Taking into account that \(|S|_{n}\leq \sqrt{n}^n|B_2^n|\) and \(\|\xi\|\geq 1\) we obtain
\[
\frac{\|\xi_{\delta}\|}{\|\xi\|}\geq 1-\sqrt{n}\left(\frac{n|B_2^n|_n}{|B_2^{n-1}|_{n-1}}\right)^{1/n}\delta^{1/n}.
\]
The desired result follows.  \hfill \(\Box\)

\vskip 3mm
\noindent
One might ask if the convergence result in Theorem \ref{MT} is uniform, i.e., does 
\[
|\d_P(\delta)-1-G(P)\delta^{1/n}|\leq o(\delta^{1/n})
\]
hold with an error term \(o(\delta^{1/n})\) only depending  on the dimension? This is not the case.
Indeed,  the floating and illumination bodies are stable with respect to the Hausdorff distance,  \cite{MeyerSchuettWerner} and hence, with respect to the distance \(\d\).
This means that  if 
\[
\lim_{n\rightarrow\infty}\d(K_n, K)=1, 
\] 
then
\[
\lim_{n\rightarrow\infty}\d\left((K_n)_{\delta},K_{\delta}\right)=1  \qquad   \text{and}  \qquad  \lim_{n\rightarrow \infty}\d\left((K_n)^{\delta},K\right)=1.
\]
Consider now the polytopes \(P(\varepsilon) \). By  continuity of the floating and illumination body we get for fixed \(\delta>0\)
\[
\lim_{\varepsilon\rightarrow 0}\frac{\d_{P(\varepsilon)}(\delta)-1}{\delta^{1/n}}=\frac{\d_{B_1^n}(\delta)-1}{\delta^{1/n}}.
\] 
On the other hand,  by  Section  \ref{Combinatorics}, 
\[
\lim_{\varepsilon\rightarrow 0}\lim_{\delta\rightarrow 0}\frac{\d_{P(\varepsilon)}(\delta)-1}{\delta^{1/n}}=\frac{2\sqrt{2}}{3}>\frac{\sqrt{2}}{2}=\lim_{\delta\rightarrow 0}\frac{\d_{B_1^n}(\delta)-1}{\delta^{1/n}}.
\]
\vskip 3mm
\noindent
We also like to address the problem to compute the optimal constant \(\tilde{G}(P)\) such that 
\[
\d_P(\delta)\leq 1+\tilde{G}(P)\delta^{1/n}+o(\delta^{1/n})
\]
for centrally symmetric polytopes such that \(o(\delta^{1/n})\) is only a dimension dependent error. The problem of proving such a result is already 
illustrated by the example  \(P(\varepsilon) \)  of Section  \ref{Combinatorics}. The facet and vertex structure of a polytope is not stable with respect  to the distance \(d\) but the convergence result Theorem \ref{MT} depends highly on these quantities. On the other hand \(P(\varepsilon) \) is close to \(B_1^n\) for small \(\varepsilon\) and therefore, 
\(d_{P(\varepsilon)}(\delta)\) and \(\d_{B_1^n}(\delta)\) behave similarly for a wide range of \(\delta\) not to close to \(0\). 
Deriving a uniform bound would demand techniques which take into account the global structure of the convex bodies.   
\section*{Acknowledgments} 
The first author would like to thank the Department of Mathematics at Case Western Reserve University, Cleveland, for their hospitality during his research stay in 2015/2016.
Both authors want to thank the Mathematical Science Research Institute, {\em Berkeley}. It was during a stay there when the paper was completed.
We also want to thank the referees for the careful reading and suggestions for improvement.

\bibliographystyle{amsplain}


\begin{dajauthors}
\begin{authorinfo}[OlMo]
  Olaf Mordhorst\\
  Goethe-Universit\"at\\
  Frankfurt am Main, Germany\\
  mordhorst\imageat{}math\imagedot{}uni-frankfurt\imagedot{}de \\
  \url{http://www.uni-frankfurt.de/76607468/ContentPage_76607468?}
\end{authorinfo}
\begin{authorinfo}[ElWe]
  Elisabeth M. Werner\\
  Case Western Reserve University\\
  Cleveland, Ohio, USA\\
  and\\
  Universit\'{e} de Lille 1\\
  Lille, France\\
  elisabeth\imagedot{}werner\imageat{}case.edu \\
  \url{https://case.edu/artsci/math/werner/}
\end{authorinfo}
\end{dajauthors}

\end{document}